\documentclass{amsart}

\theoremstyle{definition}

\theoremstyle{remark}

\newcommand{\ndiv}{\hspace{-4pt}\not|\hspace{2pt}}

\begin{document}

\title{On Fermat's Last Theorem for odd prime exponents}

\author{Zenon B. Batang}
\address{Coastal and Marine Resources Core Lab, King Abdullah University of Science and Technology, Thuwal 23955, Saudi Arabia}

\email{zenon.batang@kaust.edu.sa}


\date{April 1, 2023}


\begin{abstract}
We show that the Fermat equation $x^p + y^p = z^p$ has no solutions in coprime positive integers $x, y, z$ for any odd prime $p$.
\end{abstract}

\maketitle

For positive integers $x, y, z, n$, where $0<x<y<z$, Fermat's Last Theorem (FLT) states that
\begin{equation} \label{eq: 1}
	x^n + y^n = z^n
\end{equation}
has no positive integer solutions $(x, y, z)$ for $n>2$. For over $350$ years, FLT has been shown to hold for specific exponents \cite{Rib79} until Andrew Wiles \cite{Wil95} gave a celebrated proof in $1995$ using the machinery of modern mathematics. When Pierre de Fermat posed the problem in $1637$, he noted in a page margin of Diophantus's \textit{Arithmetica} that he possessed ``a truly remarkable proof, but the margin is too narrow to contain it." While Fermat himself proved the FLT for $n=4$ using his technique of infinite descent \cite{Bus18}, it remains the subject of debate if he indeed had a proof for the general case.

Consider the identity
\begin{equation} \label{eq: 2}
	z^n - y^n = (z-y)(z^{n-1}+y^{n-1}) + zy(z^{n-2} - y^{n-2})
\end{equation}
and let
\begin{equation} \label{eq: 3}
	z-y = a, \qquad z-x = b.
\end{equation}
It is trivial to find infinitely many solutions to \eqref{eq: 1} for $n=1$. For $n=2$, one can generate infinitely many Dickson's Pythagorean triples $(x, y, z)$ \cite{Dic94} with
\begin{equation} \label{eq: 4}
	x = a + \sqrt{2ab}, \qquad y = b + \sqrt{2ab}, \qquad z = a + b + \sqrt{2ab}.
\end{equation}
To show this, combine \eqref{eq: 1} and \eqref{eq: 2} such that
\[
\begin{aligned}
	x^2 = (z-y)(z+y) &= a(2x+2b-a) \\
	x^2-2ax+a^2 &= 2ab \\
	(x-a)^2 & = 2ab,
\end{aligned}
\]
from which \eqref{eq: 4} follows easily after substitutions into \eqref{eq: 3}. A visual demonstration of Dickson's Pythagorean triples is found in \cite{Ruc13}. Apparently, the resulting triple is primitive, i.e. $\gcd(x, y, z)=1$, if $a$ and $b$ are coprime and yield exact integer values for $\sqrt{2ab}$. For $n>2$, we will demonstrate the claim of FLT that \eqref{eq: 1} is unsolvable in positive integers by elementary number theoretic arguments.

Let
\begin{equation} \label{eq: 5}
	\phi_n(z, y) = \dfrac{z^n - y^n}{z - y} = \sum_{i=0}^{n-1} z^{n-i-1}y^i,
\end{equation}
such that
\begin{equation} \label{eq: 6}
	x^n = z^n - y^n = (z-y) \phi_n(z, y).
\end{equation}
Henceforth, we only consider the Fermat equation for odd prime $p$
\begin{equation} \label{eq: 7}
	x^p + y^p = z^p,
\end{equation}
since FLT holds for $n=4$ and for any positive integer $m$ we have that
\[
	(x^m)^p + (y^m)^p = (z^m)^p, \quad \forall  n = mp.
\]
Put $k = (p-1)/2$ for a simple manipulation of \eqref{eq: 5} into
\begin{equation} \label{eq: 8}
\begin{aligned}
	 \phi_p(z, y) &= (zy)^k + \sum_{i=0}^{k-1} (zy)^{i} (z^{2(k-i)} + y^{2(k-i)})\\
	 	&= p(zy)^k + \sum_{i=0}^{k-1} (zy)^{i} (z^{k-i} - y^{k-i})^2.
\end{aligned}
\end{equation}
For brevity, denote by
\begin{equation} \label{eq: 9}
\begin{aligned}
	 A_p(z, y) &= \sum_{i=0}^{k-1} (zy)^{i} (z^{2(k-i)} + y^{2(k-i)}),\\
	 D_p(z, y) &= \sum_{i=0}^{k-1} (zy)^{i} (z^{k-i} - y^{k-i})^2
\end{aligned}
\end{equation}
to restate \eqref{eq: 8} as
\begin{equation} \label{eq: 10}
	\phi_p(z, y) = \dfrac{z^p - y^p}{z-y} =  A_p(z, y) + (zy)^k = D_p(z, y)  + p(zy)^k.
\end{equation}
Note that \eqref{eq: 8}-\eqref{eq: 10} apply analogously to $\phi_p(z, x)$, $A_p(z, x)$ and $D_p(z, x)$. Since $p$ is odd, we also have that
\begin{equation} \label{eq: 11}
\begin{aligned}
	  \phi_p(x, -y) &= \dfrac{x^p + y^p}{x+y} =  \sum_{i=0}^{p-1} x^{p-i-1}(-y)^i \\
	  	&= A_p(x, -y) + (-1)^k(xy)^k = D_p(x, -y)  + (xy)^k,\\
\end{aligned}
\end{equation}
where we similarly obtain
\begin{equation} \label{eq: 12}
\begin{aligned}
	 A_p(x, -y) &= \sum_{i=0}^{k-1} (-1)^i (xy)^{i} (x^{2(k-i)} + y^{2(k-i)}), \\
	 D_p(x, -y) &= \sum_{i=0}^{k-1} (-1)^i (xy)^{i} (x^{k-i} - y^{k-i})^2.
\end{aligned}
\end{equation}

Given $xyz \neq 0$, $\gcd(x, y, z) = 1$ and $\gcd(p, xyz) = 1$ in \eqref{eq: 7}, it is known \cite{Bir04} that
\begin{equation} \label{eq: 13}
	\gcd(z-y, \phi_p(z,y)) =1.
\end{equation}
To prove this, assume for contradiction that there exists a prime $q$ dividing both $(z-y)$ and $\phi_p(z,y)$. Substituting $z \equiv y\pmod{q}$ into \eqref{eq: 8} using $A_p(z, y)$ or $D_p(z, y)$ yields $\phi_p(z, y) \equiv py^{p-1} \pmod{q}$. As $q \vert \phi_p(z, y)$, we have either $q=p$ or $q \vert y$. If $q=p$, then $p \vert (z-y)$ by assumption and, hence, $p \vert x$, which contradicts $\gcd(p, xyz)=1$. If $q \vert y$, then $q$, different from $p$, is a common factor of $y$ and $z$, thus also contradicting $\gcd(x, y, z) = 1$. Therefore, both cases imply that \eqref{eq: 13} must hold.

We also claim that
\begin{equation} \label{eq: 14}
	z-y \equiv 0 \pmod{p} \quad \Longrightarrow \quad \phi_p(z,y) \equiv 0 \pmod{p},
\end{equation}
since if $p \vert (z-y)$, then $p \vert D_p(z,y)$ as $(z-y)$ is a common factor of all summands in $D_p(z,y)$. With $p$ being a factor of $p(zy)^k$, then $p \vert \phi_p(z,y)$, as claimed. It follows that if $p \vert (z-y)$ or $p \vert (z-x)$, then by \eqref{eq: 14} it is assured that $p^2 \vert x^p$ or  $p^2 \vert y^p$, respectively. Moreover, given \eqref{eq: 13}, it is apparent that $\gcd(z-y, \phi_p(z,y)) = 1$ or $p$.

Following the approach of Sophie Germain \cite{Lau10}, we split FLT into two cases, where the Fermat equation has no nontrivial integer solutions for which $p \ndiv xyz$ (FLT1) or $p \vert xyz$ (FLT2). By Fermat's Little Theorem, any positive integer $N$ that is coprime to $p$ satisfies
\[
	N^p \equiv N \pmod{p} \quad \Longrightarrow \quad N^{p-1} \equiv 1 \pmod{p}.
\]
If FLT1 fails, such that Fermat equation has a solution for $p$ under FLT1 conditions, i.e. $\gcd(x, y, z) = 1$ and $p \ndiv xyz$, then
\begin{equation} \label{eq: 15}
	z^p -y^p \equiv x^p \pmod{p} \quad \Longrightarrow \quad z-y \equiv x \pmod{p}.
\end{equation}
We know from \eqref{eq: 13} that $(z-y)$ and $\phi_p(z, y)$ are coprime, such that for some positive integers $r, r_1$, where $p \ndiv r r_1$, and analogously for some positive integers $s, s_1$ and $t, t_1$, where $p \ndiv s s_1$ and $p \ndiv t t_1$, one can write the Barlow-Abel relations \cite{Rib79} as
\[
\begin{aligned}
	z-y &= r^p, \qquad && \phi_p(z, y) = r_1^p, \qquad && x = r r_1, \\
	z-x &= s^p, \qquad && \phi_p(z, x) = s_1^p, \qquad && y = s s_1, \\
	x+y &= t^p, \qquad &&\phi_p (x, -y) = t_1^p, \qquad && z = t t_1.
\end{aligned}
\]
If along with \eqref{eq: 3} we let $x+y = c$, then $a,b,c$ are all nonzero $p$-th power integers if Fermat equation holds for $p$ under FLT1 conditions.

Now we extend Dickson's method for $n > 1$ with the formula
\begin{equation} \label{eq: 16}
	(x-a)^n = f_n(x,b,a) + g_n(b,a),
\end{equation}
with
\begin{equation} \label{eq: 17}
\begin{aligned}
	f_n(x,b,a) & = - \sum_{i=1}^{n-1} \binom{n}{i}b^i \bigg((x-a)^{n-i} - (x^{n-i} + (-a)^{n-i})\bigg),\\
	g_n(b,a) &= - (b-a)^n + b^n + (-a)^n = - \sum_{i=1}^{n-1} \binom{n}{i} b^{n-i}a^i,
\end{aligned}
\end{equation}
where $\binom{n}{i}$ denotes the binomial coefficient and $a, b$ are as defined in \eqref{eq: 3}. Note that the binomial coefficients appear in all resulting terms of $f_n(x,b,a)$ and $g_n(b,a)$. For some integers $u, v$, where $\gcd(p, uv)=1$, it always holds that
\[
	u^p \pm v^p \equiv  (u \pm v)^p \pmod{p},
\]
since, by the Binomial Theorem, we have the expansion
\[
	(u \pm v)^ p = u^p \pm \binom{p}{1} u^{p-1} v + \binom{p}{2}  u^{p-2} v^2 \pm \cdots + \binom{p}{p-1} u v^{p-1} \pm v^p,
\]
where only the first and last terms are not divisible by $p$. Hence, we can claim that
\begin{equation} \label{eq: 18}
	(x-a)^p \equiv 0 \pmod{p} \qquad \Longrightarrow \qquad x \equiv a \pmod{p},
\end{equation}
which coincides with \eqref{eq: 15} if FLT1 fails, where $a=z-y$. From \eqref{eq: 3}, we obtain
\begin{equation} \label{eq: 19}
	x - a = x + y - z , \qquad 2x + b - a = x + y, \qquad y = x + b - a,
\end{equation}
such that any inferences on $x$ from \eqref{eq: 16} apply analogously to $y$ and $z$ since
\begin{equation} \label{eq: 20}
	x - a = y - b = z - (a + b) = x + y - z.
\end{equation}
Therefore, we can express
\begin{equation} \label{eq: 21}
	(x-a)^p = (y-b)^p = \big(z-(a+b)\big)^p = pabK_p,
\end{equation}
where
\begin{equation} \label{eq: 22}
\begin{split}
	K_p &= - \sum_{i=1}^{p-1} \dfrac{1}{i} \binom{p-1}{i-1}b^{i-1} \bigg( \sum_{j=1}^{p-i-1}(-1)^j \binom{p-i}{j}x^{p-i-j} a^{j-1} \bigg) \\
			& \qquad \qquad \quad -  \sum_{i=1}^{p-1} \dfrac{(-1)^i}{i} \binom{p-1}{i-1}b^{p-i-1}a^{i-1}.
\end{split}
\end{equation}
Notice that $pab \vert (x-a)^p$ but $(pab)^2 \ndiv (x-a)^p$ for all $p$, since $pab$ is the only common factor of the resulting terms after expansion of $f_p(x,b,a)$ and $g_p(b,a)$. Hence, it is apparent that
\begin{equation} \label{eq: 23}
	(x-a)^p = (y-b)^p = \big(z-(a+b)\big)^p \equiv 0 \pmod{p}
\end{equation}
and so
\begin{equation} \label{eq: 24}
	x \equiv a, \qquad y \equiv b, \qquad z \equiv a+b \pmod{p}.
\end{equation}

For the Fermat equation to have a solution for $p$ under FLT1 conditions, it must hold that $p \ndiv a$ or $p \ndiv b$ for, otherwise, we have $p \vert (z-y)$ or $p \vert (z-x)$, which contradicts our assumption that $p \ndiv xyz$. From \eqref{eq: 21}, we have that
\begin{equation} \label{eq: 25}
	x - a = y - b = z - (a+b) = \sqrt[\leftroot{-4} \uproot{5} p]{pabK_p} = rs \sqrt[\leftroot{-4} \uproot{5} p]{pK_p},
\end{equation}
where $r, s$ are as defined for the first Barlow-Abel relations above. This shows that the $p$-th root of $pK_p$ cannot be an integer since $p \ndiv K_p$. Observe that
\[
	g_p(b,a) = - pab \sum_{i=1}^{p-1} \dfrac{(-1)^i}{i} \binom{p-1}{i-1}b^{p-i-1}a^{i-1},
\]
such that
\[
	(x-a)^p = f_p(x,b,a) + g_p(b,a) \qquad \textnormal{and} \qquad K_p = \dfrac{f_p(x,b,a)}{pab} +  \dfrac{g_p(b,a)}{pab}.
\]
Then substituting $x \equiv a \pmod{p}$ into \eqref{eq: 17} yields
\begin{equation} \label{eq: 26}
\begin{split}
	K_p & \equiv 2 \sum_{i=1}^{k} \dfrac{1}{2i-1} \binom{2k}{2i-2}b^{2(i-1)} a^{2(k-i)+1} \\
			& \qquad - \sum_{j=1}^{2k} \dfrac{(-1)^j}{j} \binom{2k}{j-1}b^{2k-j}a^{j-1} \pmod{p},
\end{split}
\end{equation}
where, as above, $k = (p-1)/2$. As the exponents are zero when $i=1$ and $j=2k$, the right-hand side of \eqref{eq: 26} forms a polynomial with all terms consisting of different powers of $a$, $b$, or $ab$. This implies $p \ndiv K_p$ from \eqref{eq: 25}, with the radical factor assuring that $(x-a)$ is not an integer, thus establishing the truth of FLT1 for all $p$.

On the other hand, suppose FLT2 fails, such that Fermat equation has a solution for $p$ under FLT2 conditions, i.e. $\gcd(x, y, z) = 1$ and $p \vert xyz$. Since one of $x, y, z$ is necessarily even, we consider two cases where $p$ divides the even (Case I) or one of the odd (Case II) variables. Without loss of generality, assume that $x$ is even and divisible by $p$, such that $2p \vert a$ or $2p \vert (z-y)$. Then we have that
\begin{equation} \label{eq: 27}
	z^p - y^ p \equiv 0 \pmod{p} \qquad \Longrightarrow \qquad z \equiv y \pmod{p}.
\end{equation}
From \eqref{eq: 14}, it follows that $\phi_p(z,y) \equiv 0 \pmod{p}$. Since $x^p = (z-y)\phi_p(z,y)$ and, from \eqref{eq: 8}, $p \vert \phi_p(z,y)$ but $p^2 \ndiv \phi_p(z,y)$, the Barlow-Abel relations \cite{Rib79} for Case I must be in the form 
\[
\begin{aligned}
	z-y &= 2^{pd}p^{pe-1}r'^p, \qquad && \phi_p(z, y) = pr_1^p, \qquad && x = 2^{d}p^{e}r' r_1, \\
	z-x &= s^p, \qquad && \phi_p(z, x) = s_1^p, \qquad && y = s s_1, \\
	x+y &= t^p,  \qquad && \phi_p (x, -y) = t_1^p, \qquad && z = t t_1,
\end{aligned}
\]
for some positive integers $d, e, r'$, where $2p \ndiv r'$. Similarly for Case II, assuming that $p$ divides $y$ instead of $x$, we have that
\[
\begin{aligned}
	z-y &= 2^{pd} r_0^p, \qquad && \phi_p(z, y) =  r_1^p, \qquad && x = 2^{d} r_0 r_1, \\
	z-x &= p^{pe-1} s'^p, \qquad && \phi_p(z, x) = p s_1^p, \qquad && y = p^e  s' s_1, \\
	x+y &= t^p,  \qquad && \phi_p (x, -y) = t_1^p, \qquad && z = t t_1,
\end{aligned}
\]
for some positive integers $d, e, r_0, s'$, where $2p \ndiv r_0s'$. One can see that $a=z-y$ and $b = z-x$ are not perfect $p$-th powers for Case I and Case II, respectively. Hence, FLT2 conditions also result in \eqref{eq: 21}, but where
\begin{equation} \label{eq: 28}
	\sqrt[\leftroot{-4} \uproot{5} p]{pabK_p} =
\begin{cases}
	2^d p^e r' s \sqrt[\leftroot{-4} \uproot{5} p]{K_p}, \quad & \textnormal{for FLT2 Case I }, \\
	2^d p^e r_0 s' \sqrt[\leftroot{-4} \uproot{5} p]{K_p}, \quad & \textnormal{for FLT2 Case II}.
\end{cases}
\end{equation}
Since $p \vert a$ in Case I or $p \vert b$ in Case II, we can rewrite \eqref{eq: 26} as
\begin{equation} \label{eq: 29}
	K_p \equiv h_p(b,a) + a^{p-2} + b^{p-2}  \pmod{p},
\end{equation}
where it can be checked that $h_p(b,a)=0$ for $p=3$ and $a \vert h_p(b,a)$ for $p>3$ and so
\begin{equation} \label{eq: 30}
	K_p \equiv 
\begin{cases}
	b  & \textnormal{for FLT2 Case I},\\
	a  & \textnormal{for FLT2 Case II},
\end{cases}
	\pmod{p},
\end{equation}
which implies that $p \ndiv K_p$. We argue that $K_p$ is never a perfect $p$-th power, except for the special case when $K_p=1$ for $p=2$ or also known as Pythagorean Theorem.

We can show from the Barlow-Abel relations above that
\begin{equation} \label{eq: 31}
	2x = c - b + a, \qquad 2y = c + b - a, \qquad 2z = c + b + a,
\end{equation}
which also lead to \eqref{eq: 20}. Obviously, it is guaranteed that
\[
	x - a = \sqrt[\leftroot{-4} \uproot{5} p]{pabK_p}
\]
holds when $K_p = 1$ for $p=2$, thus making the Pythagorean Theorem an exception to FLT. This is indeed a special case, as only Case I of FLT2 is valid to apply for $p=2$. From \eqref{eq: 16}, we obtain $f_2(x,b,a) = 0$ and easily recover
\begin{equation} \label{eq: 32}
	(x-a)^2 = -(b-a)^2 + b^2+a^2 = 2ab,
\end{equation}
which leads to Dickson's Pythagorean triple in \eqref{eq: 4}, since from \eqref{eq: 20} we also have that $(y-b)^2= 2ab$ and $(z-(a+b))^2 = 2ab$. Then, from \eqref{eq: 21}, we get $(x-a)^2 = 2abK_2$, where $K_2 = 1$. The Barlow-Abel relations for Case I of FLT2 must take $a = 2^{2e-1}r'^2$ and $b = s^2$, where $2 \ndiv r's$, such that from \eqref{eq: 28} we obtain
\[
	x-a = \sqrt{2abK_2} = 2^{e} r's,
\]
whence Dickson's Pythagorean triple follows accordingly. This applies analogously if one assumes that $y$, instead of $x$, is even.

On the assumption that Fermat equation holds for odd $p$, it is claimed that
\begin{equation} \label{eq: 33}
	x^p + y^p = z^p \qquad \Longrightarrow \qquad  (x-a)^p - pabK_p = 0.
\end{equation}
Hence, \eqref{eq: 16} is just a reformulation of Fermat equation with the premise that
\begin{equation} \label{eq: 34}
	(x+y-z)^p - pabK_p = x^p + y^p - z^p = 0.
\end{equation}
In other words, $pabK_p$ is the summation of all multivariate terms after expansion of $(x+y-z)^p$ for odd $p$, with the necessary condition that $pabK_p$ is a perfect $p$-th power for the Fermat equation to be true. However, FLT asserts that $x^p + y^p -z^p \ne 0$ in positive integers for odd $p$ and the Barlow-Abel relations above were, in fact, based on a contrary assumption that FLT fails. With the eventual proof of FLT by Wiles \cite{Wil95}, it thus hints at rethinking the admissibility of the Barlow-Abel relations above due to a false assumption, which we illustrate in our argument below.

For $p=3$, we find that $f_3(x,b,a) = 6abx$ and $g_3(b,a) = -3ab^2+3a^2b$, such that
\begin{equation} \label{eq: 35}
	(x-a)^3 =  6abx + 3ab^2 - 3a^2b = 3ab(2x+b-a).
\end{equation}
One can check from \eqref{eq: 22} that $K_3 = 2x+b-a$ and so, given \eqref{eq: 24}, we get $K_3 \equiv a+b \pmod{3}$ and, hence, $3 \ndiv K_3$ since $3 \ndiv ab$ under FLT1 conditions. It should hold that $3^3 \vert (x-a)^3$ and, by unique factorization, $3^2 \vert ab(x+y)$ where $K_3=x+y$ from \eqref{eq: 19}. Given that $3 \ndiv ab$, with $a$ and $b$ being perfect cubes by the Barlow-Abel relations, it follows that $3^2 \vert (x+y)$. But then this implies that $z \equiv 0 \pmod{3}$, which is contrary to our assumption that $3 \ndiv xyz$. Hence, this tells us that FLT1 must hold for $p=3$, as there are no integer solutions to
\[
	x-a = rs \sqrt[\leftroot{-4} \uproot{3} 3]{3K_3}.
\] 
On the other hand, both cases of the Barlow-Abel relations under FLT2 conditions lead to
\[
	x-a = y-b = z -(a+b) = \sqrt[\leftroot{-4} \uproot{5} 3]{3abK_3},
\]
where, from \eqref{eq: 28}, we have that
\begin{equation} \label{eq: 36}
	\sqrt[\leftroot{-4} \uproot{5} 3]{3abK_3} =
\begin{cases}
	2^d 3^e r' s \sqrt[\leftroot{-2} \uproot{3} 3]{K_3}, \quad & \textnormal{for FLT2 Case I}, \\
	2^d 3^e r_0 s' \sqrt[\leftroot{-2} \uproot{3} 3]{K_3}, \quad & \textnormal{for FLT2 Case II}.
\end{cases}
\end{equation}
Given that $z=x+b$, observe that
\[
	z = K_3 - (x-a),
\]
which is only valid if $(x-a)$ is an integer, for which $K_3$ must be a perfect cube. We argue that the claim in the Barlow-Abel relations, where $x+y = t^p$, is inadmissible if one makes the stated assumptions for $a=z-y$ and $b=z-x$. Indeed, it can be shown that
\[
	(x+y-z)^3 = 3(z-y)(z-x)(x+y),
\]
but a simple algebraic manipulation leads us to deduce that
\[
	\dfrac{(x+y-z)^3}{3(z-y)(z-x)} \ne t^3.
\]
This suggests that $K_3=x+y$ is not a cube, such that $(x-a)$ cannot be an integer if one assumes the Fermat equation. Hence, FLT2 must hold for $p=3$.

For $p=5$, we find that
\begin{equation} \label{eq: 37}
\begin{split}
	f_5(x,b,a) &= 5ab (4x^3 - 6ax^2 + 4a^2x + 6bx^2 + 4b^2x - 6abx), \\
	g_5(b,a) &= 5ab (b^3 - 2ab^2 + 2a^2b - a^3),
\end{split}
\end{equation}
from which we can formulate
\begin{equation} \label{eq: 38}
	 K_5 = 4x^3 + 6(b-a)x^2 + 2(2b^2 -ab + 2a^2)x - 2abc + b^3 - a^3,
\end{equation}
with all terms being even, except for $b^3$, and thus $K_5$ is odd. By substituting $x \equiv a \pmod{5}$ into $f_5(x,b,a)$, we obtain 
\begin{equation} \label{eq: 39}
	K_5 \equiv h_5(b,a) + b^3 + a^3 \pmod{5},
\end{equation}
where one can check that $h_5(b,a) = 2ab^2 + 2a^2b$. This implies that $5 \ndiv K_5$ and FLT1 must hold for $p=5$ since there are no integer solutions to
\[
	x-a = y-b = z -(a+b)= \sqrt[\leftroot{-4} \uproot{5} 5]{5abK_5} = rs \sqrt[\leftroot{-4} \uproot{5} 5]{5K_5}.
\]
Under FLT2 conditions, suppose $K_p$ is a perfect $p$-th power, such that there exists an integer $R$ satisfying $R^p = pabK_p$. Then $x^p = (a+R)^p$ and $a^p = (z-y)^p$ lead to
\begin{equation} \label{eq: 40}
	x^p + y^p - z^p = R^p + \sum_{i=1}^{p-1} \binom{p}{i} \bigg( (z-y)^{p-i}R^i + (-z)^{p-i} y^i \bigg),
\end{equation}
which is clearly nonzero and, hence, in contradiction to \eqref{eq: 34}. Therefore, FLT2 must also hold for $p=5$ since it is impossible to find integer solutions to
\[
	\sqrt[\leftroot{-4} \uproot{5} 5]{5abK_5} =
\begin{cases}
	2^d 5^e r' s \sqrt[\leftroot{-2} \uproot{3} 5]{K_5}, \quad & \textnormal{for FLT2 Case I}, \\
	2^d 5^e r_0 s' \sqrt[\leftroot{-2} \uproot{3} 5]{K_5}, \quad & \textnormal{for FLT2 Case II}.
\end{cases}
\]
By induction on $K_p$ using \eqref{eq: 40}, we thus conclude that FLT is true for all odd $p$.

\bibliographystyle{amsplain}

\end{document}